\RequirePackage{ifpdf}
\ifpdf % We are running pdfTeX in pdf mode
\documentclass[pdftex]{sigma}
\else
\documentclass{sigma}
\fi

\begin{document}

\allowdisplaybreaks

\renewcommand{\thefootnote}{$\star$}

\renewcommand{\PaperNumber}{025}

\FirstPageHeading

\ShortArticleName{Inversion Formulas for the Spherical Radon--Dunkl
Transform}

\ArticleName{Inversion Formulas\\ for the Spherical Radon--Dunkl
Transform\footnote{This paper is a contribution to the Special
Issue on Dunkl Operators and Related Topics. The full collection
is available at
\href{http://www.emis.de/journals/SIGMA/Dunkl_operators.html}{http://www.emis.de/journals/SIGMA/Dunkl\_{}operators.html}}}

\Author{Zhongkai LI and Futao SONG}

\AuthorNameForHeading{Zh.-K. Li and F.-T. Song}

\Address{Department of Mathematics, Capital Normal University,
Beijing 100048, China}
\Email{\href{mailto:lizk@mail.cnu.edu.cn}{lizk@mail.cnu.edu.cn}, \href{mailto:ftsong@tom.com}{ftsong@tom.com}}

\ArticleDates{Received October 18, 2008, in f\/inal form March 01,
2009; Published online March 03, 2009}

\Abstract{The spherical Radon--Dunkl transform $R_{\kappa}$,
associated to weight
 functions invariant under a f\/inite ref\/lection group, is introduced,
 and some elementary properties are obtained in terms of
$h$-harmonics. Several inversion formulas of $R_{\kappa}$ are given
with the aid of spherical Riesz--Dunkl potentials, the Dunkl
operators, and some appropriate wavelet transforms.}


\Classification{44A12; 33C55; 65R32; 42C40}

%%%%%%%%%%%%%%%%%%%%%%%%%%%%%%%%%%%%%%%%%%%%%%%%%%%%%%%%%%%%%%%%%%%%%%%%%%%%%%%%%%%%%%%%%%%%%

%%%%%%%%%%%%%%%%%%%%%%%%%%%%%%%%%%%%%%%%%%%%%%%%%%%%%%%%%%%%%%%%%%%%%%%%%%%%%%%%%%%%%%%%%%%%%
\section{Introduction}

Let $\langle x,y\rangle$ denote the usual Euclidean inner product of
$x,y\in{\mathbb R}^{d+1}$, and ${\mathbb S}^d=\{x: \; \|x\|=1\}$ the
unit sphere in ${\mathbb R}^{d+1}$. We use $d\omega_{k}$ to denote
the surface (Lebesgue) measure on a $k$-dimensional sphere. The
spherical Radon transform $R$ is one of the tools in integral
geometry, which is def\/ined, for $f\in C({\mathbb S}^d)$, by
\begin{gather*}
Rf(x)=\Lambda_{d-1}\int_{\langle x,y\rangle=0}f(y)d\omega_{d-1}(y),
\qquad x\in{\mathbb S}^d,
\end{gather*}
where $\Lambda_{d-1}^{-1}$ is the surface area of ${\mathbb
S}^{d-1}$. There are a number of papers devoting to the study of the
spherical Radon transform $R$ by dif\/ferent methods (see \cite{Ca,Gr,Gu,He1,He2,Ku,Ru2,Ru3,Ru4,Ru5,Ru6,Ru7,RR,St2}), and to its
applications to various problems (see \cite{GRS,GW}). Some
deep results about $R$ were obtained with the aid of spherical
harmonics (see \cite{Ru2,Ru3,Ru4,Ru5,Ru6,Ru7,RR,St2}), and furthermore,
$R$ is a special case of the spherical means
\begin{gather*}
M_{\tau}f(x)=\frac{\Lambda_{d-1}}{(1-\tau^2)^{(d-1)/2}}\int_{\langle
x,y\rangle=\tau}f(y)d\omega_{d-1}(y), \qquad x\in{\mathbb S}^d,
\end{gather*}
by taking $\tau=0$, and also of the spherical Riesz potentials
\begin{eqnarray*}
I^{\alpha}f(x)=\frac{\Gamma((1-\alpha)/2)}{2\pi^{d/2}\Gamma(\alpha/2)}\int_{{\mathbb
S}^d} |\langle x,y\rangle|^{\alpha-1}f(y)d\omega_{d}(y), \qquad
x\in{\mathbb S}^d,
\end{eqnarray*}
by taking the limit $\lim\limits_{\alpha\rightarrow0+}I^{\alpha}f=Rf$ in
some sense (see~\cite{Ru3}). The former is a tool in approximation
on the sphere ${\mathbb S}^d$, and the later is one of the research
objectives in harmonic analysis on~${\mathbb S}^d$.

The purpose of the present paper is to study an analogous model of
the spherical Radon transform $R$ in Dunkl's theory. This is based
on the def\/inition of the generalized spherical means
$M_{\tau}^{\kappa}f(x)$ due to~\cite{Xu3} (instead of
$M_{\tau}^{\kappa}$, the notation $T_{\theta}^{\kappa}$ with
$\tau=\cos\theta$ was used there), in terms of the equation
\begin{gather*}
\int_{-1}^1
M_{\tau}^{\kappa}f(x)g(\tau)w_{\lambda_{\kappa}}(\tau)d\tau
=c_{\kappa}\int_{{\mathbb S}^d}f(y)V_{\kappa}[g(\langle
x,\cdot\rangle)](y)h_{\kappa}^2(y)d\omega_d(y)
\end{gather*}
for any $g$ in $L^1([-1,1];w_{\lambda_{\kappa}})$, where
$w_{\lambda_{\kappa}}(t)=
\tilde{c}_{\lambda_{\kappa}+1/2}(1-t^2)^{\lambda_{\kappa}-1/2}$,
$V_{\kappa}$ is the intertwining operator associated to a given
f\/inite ref\/lection group, and $h_{\kappa}^2$ is the related weight
function (for details concerning them and other notations in the
equation, see the next section). We def\/ine~$R_{\kappa}$ by
$R_{\kappa}f=M_0^{\kappa}f$, and call $R_{\kappa}$ the spherical
Radon--Dunkl transform. Although~$M_{\tau}^{\kappa}$ is def\/ined
implicitly, it is a proper extension of $M_{\tau}$ and
$M_{\tau}^0=M_{\tau}$, and moreover, from~\cite{DX} and
\cite{Xu3,Xu4,Xu5,Xu6}, $M_{\tau}^{\kappa}$ shares many properties
with $M_{\tau}$ and plays the same roles in weighted approximation
and related harmonic analysis on the sphere ${\mathbb S}^d$. One
could expect that the spherical Radon--Dunkl Transform $R_{\kappa}$
would have similar features to $R$ and be a suitable tool in
reconstruction of functions in weighted spaces. This is the
motivation of the paper. Despite less closed representation, a
further work worth doing is to f\/ind applications of $R_{\kappa}$ in
geometry or other f\/ields.

The paper is organized as follows.
In Section~\ref{section2}, some necessary facts in Dunkl's theory are reviewed,
and in Section~\ref{section3}, the spherical Radon--Dunkl transform $R_{\kappa}$
is def\/ined and some of elementary properties are obtained in terms
of $h$-harmonics. Sections~\ref{section4} and~\ref{section5} are devoted to inversion formulas
of~$R_{\kappa}$, which are given by means of spherical Riesz--Dunkl
potentials~$I_{\kappa}^{\alpha}$, the Dunkl operators, and some
appropriate wavelet transforms. These conclusions generalize part of
those in~\cite{Ru2,Ru3,Ru4}.

\section{Some facts in Dunkl's theory}\label{section2}

Let $G$ be a f\/inite ref\/lection group on ${\mathbb R}^{d+1}$ with a
f\/ixed positive root system $R_+$, normalized so that $\langle
v,v\rangle=2$ for all $v\in R_+$. It is known that $G$ is a subgroup
of $O(d+1)$ generated by $\{\sigma_v: \; v\in R_+\}$, where
$\sigma_v$ denotes the ref\/lection with respect to the hyperplane
perpendicular to~$v$, i.e.\ $x\sigma_v=x-2(\langle x,v\rangle/\langle
v,v\rangle)v$ for $x\in{\mathbb R}^{d+1}$. Let $\kappa$ be a
multiplicity function $v\mapsto\kappa_v\in[0,+\infty)$ def\/ined on
$R_+$, with invariance under the action of $G$. Thus $\{\kappa_v: \;
v\in R_+\}$ has dif\/ferent values only as many as the number of
$G$-orbits in $R_+$.

The Dunkl operators are a family of f\/irst-order
dif\/ferential-ref\/lection operators ${\mathcal{D}}_j, \; 1\le j\le
d+1$, def\/ined by (see \cite{Du2})
\begin{gather*}
  {\mathcal{D}}_jf(x):=\partial_jf(x)+\sum_{v\in R_+}\kappa_v\frac{f(x)-
  f(x\sigma_v)}{\langle x,v\rangle}\langle v,e_j\rangle,
\end{gather*}
for $f\in C^{1}({\mathbb{R}}^{d+1}),$ where $\{e_i: \; 1\le i\le
d+1\}$ is the usual standard basis of ${\mathbb R}^{d+1}$. As
substitutes of partial dif\/ferentiations $\partial_j$, these
operators are mutually commutative. The associated Laplacian, called
$h$-Laplacian, is def\/ined by
$\Delta_h={\mathcal{D}}_1^2+\cdots+{\mathcal{D}}_{d+1}^2$, which
plays roles similar to that of the usual Laplacian $\Delta=\Delta_0$
(see \cite{Du1}). In terms of the polarspherical coordinates $x=rx'$,
$r=\|x\|$, the operator $\Delta_h$ can be expressed as (see
\cite{Xu3})
\begin{gather*}
\Delta_h=\frac{\partial^2}{\partial
r^2}+\frac{2\lambda_{\kappa}+1}{r}\frac{\partial}{\partial
r}+\frac{1}{r^2}\Delta_{h,0},
\end{gather*}
where $\Delta_{h,0}$ is the associated Laplace-Beltrami operator on
${\mathbb S}^d$, and $\lambda_{\kappa}=\gamma_{\kappa}+(d-1)/2$ with
$\gamma_{\kappa}=\sum\limits_{v\in R_+}\kappa_v$. If $\gamma_{\kappa}=0$,
i.e. $\kappa_v\equiv0$, then ${\mathcal{D}}_j=\partial_j$, $1\le
j\le d+1$. In the following, we assume that $\gamma_{\kappa}>0$, and
so $\lambda_{\kappa}>0$.

For each multiplicity function $\kappa$, there is a linear operator
$V_{\kappa}$ intertwining the partial dif\/fe\-ren\-tia\-tions and the Dunkl
operators (see~\cite{Du3}). Precisely, if
${\mathcal{P}}_n={\mathcal{P}}_n^{d+1}$ denotes the set of
homogeneous polynomials of degree $n$ in $d+1$ variables, then the
intertwining operator $V_{\kappa}$ is determined
 uniquely by
 $V_{\kappa}{\mathcal{P}}_n\subseteq {\mathcal{P}}_n$,
$V_{\kappa}1=1$ and ${\mathcal{D}}_jV_{\kappa}=V_{\kappa}\partial_j
$, $1\le i\le d+1$.
 $V_{\kappa}$ commutes with the group action and is
a linear isomorphism on each ${\mathcal{P}}_n$ $(n=0,1,\dots)$.
Moreover it is a~positive operator (see~\cite{Ro1}) and can be
extended to the space of smooth functions and even to the space of
distributions (see~\cite{Tr1,Tr2}).

The intertwining operator $V_{\kappa}$ allows to introduce some
useful tools in Dunkl's theory. For example, the Dunkl transform
${\mathcal{\mathcal{F}}}_{\kappa}$ is def\/ined in \cite{Du4}
associated with the measure $h_{\kappa}^2dx$ on ${\mathbb R}^{d+1}$,
where
\begin{gather*}
h_{\kappa}(x)=\prod_{v\in R_+}|\langle x,v\rangle|^{\kappa_v}.
\end{gather*}
${\mathcal{F}}_{\kappa}$ is a generalization of the Fourier
transform ${\mathcal{F}}={\mathcal {F}}_0$ and enjoys properties
similar to those of~${\mathcal{F}}$ (see~\cite{dJ,Du4,Ro2}).

For $1\le p<\infty$, denote by
$\|f\|_{\kappa,p}=\left\{c_{\kappa}\int_{{\mathbb S}^d}|f|^p
h_{\kappa}^2d\omega_d\right\}^{1/p}$ the norm of $f\in L^p({\mathbb
S}^d;h_{\kappa}^2)$, with $c_{\kappa}^{-1}=\int_{{\mathbb
S}^d}h_{\kappa}^2d\omega_d$, and by
$\|\phi\|_{\lambda_{\kappa},p}=\left\{\int_{-1}^1|\phi|^pw_{\lambda_{\kappa}}dt\right\}^{1/p}$
the norm of $\phi\in L^p([-1,1],w_{\lambda_{\kappa}})$, where
$w_{\lambda_{\kappa}}(t)=
\tilde{c}_{\lambda_{\kappa}+1/2}(1-t^2)^{\lambda_{\kappa}-1/2}$,
$\tilde{c}_{\lambda}=\pi^{-1/2}\Gamma(\lambda+1/2)/\Gamma(\lambda)$.
When $p=\infty$, $\|f\|_{\infty}=\|f\|_{\kappa,\infty}$ and
$\|\phi\|_{\infty}=\|\phi\|_{\kappa,\infty}$ are def\/ined as usual.

 The functions in
${\mathcal{H}}_{n}^{h,d+1}:= {\mathcal{P}}_{n}^{d+1}\cap\,
\ker\Delta_{h}$ are called $h$-harmonic polynomials of degree $n$,
and the spherical $h$-harmonics of degree $n$ are their restrictions
on ${\mathbb S}^d$. The orthogonality theorem in \cite{Du1} asserts
that if $P\in{\mathcal{P}}_n^{d+1}$, then
$\int_{{\mathbb{S}}^d}PQh_{\kappa}^{2}d\omega_d=0$ for all
$Q\in\mathop{\cup}\limits_{k=0}^{n-1}{\mathcal{P}}_{k}^{d+1}$, if and only if $P$
is $h$-harmonic, i.e.\ $\Delta_{h} P=0$. Moreover $L^2({\mathbb
S}^{d};h_\kappa^2)=\sum\limits_{n=0}^{\infty}\bigoplus{\mathcal{H}}_n^{h,d+1}$.

If $Y_n(h_\kappa^2;f;x)$ is the projection of $f\in L^1({\mathbb
S}^{d};h_\kappa^2)$ to ${\mathcal{H}}_{n}^{h,d+1}$, then the
$h$-harmonic expansion of $f$ is given by
\begin{gather}\label{eq1}
f(x)\sim\sum_{n=0}^\infty Y_n\big(h_\kappa^2;f;x\big), \qquad
x\in{\mathbb S}^d.
\end{gather}
The projection $Y_n(h_\kappa^2;f;x)$ takes the form
\begin{gather}\label{eq2}
Y_n(h_\kappa^2;f;x) = c_{\kappa}\int_{{\mathbb{S}}^{d}} f(y)
P_n\big(h_\kappa^2;x,y\big)h_\kappa^2(y) d\omega_d(y),
\end{gather}
where $P_n(h_\kappa^2;x,y)$ is the reproducing kernel of the space
${\mathcal{H}}_{n}^{h,d+1}\!$.\ A compact formula of
$P_n(h_\kappa^2;x,y)\!$ is (see~\cite{Xu2})
\begin{gather}\label{eq3}
  P_n(h_\kappa^2;x,y)=\frac{n+\lambda_{\kappa}}
    {\lambda_{\kappa}}
   V_{\kappa}\big[C_n^{\lambda_{\kappa}} (\langle x,\cdot\rangle )\big](y),
\end{gather}
with $C_n^{\lambda_{\kappa}}$, the Gegenbauer polynomial of degree
$n$ with parameter $\lambda_{\kappa}$. It is noted that (see~\cite{Xu3})
\begin{gather}\label{eq4}
\Delta_{h, 0}Y_n=-n(n+2{\lambda_\kappa})Y_n, \qquad
Y_n\in{\mathcal{H}}_{n}^{h,d+1}.
\end{gather}

When $\kappa_v=0$ for all $v\in R_+$, we have $V_0=id$, and hence,
$P_n(h_\kappa^2;x,y)$ reduces to the usual zonal polynomial for the
ordinary spherical harmonics $P_n(x,y) = \frac{n+(d-1)/2}{(d-1)/2}
C_n^{(d-1)/2}(\langle x,y \rangle )$.

A useful integration formula for the intertwining operator
$V_{\kappa}$  is
\begin{gather}\label{eq5}
\int_{{\mathbb{S}}^{d}} V_{\kappa}f(x)h_\kappa^2(x)
d\omega_d(x)=\frac{c_{\kappa}^{-1}\Gamma(\lambda_{\kappa}+1)}{\pi^{(d+1)/2}\Gamma(\gamma_{\kappa})}
\int_{{\mathbb{B}}^{d+1}} f(x)\big(1-|x|^2\big)^{\gamma_{\kappa}-1}dx.
\end{gather}
The formula is proved in \cite{Xu1} when $f$ is a polynomial.
Applying density of polynomials and positivity of $V_{\kappa}$, this
allows us to extend the intertwining operator $V_{\kappa}$ acting on
those functions~$f$ on the sphere ${\mathbb S}^d$ which are
restrictions of functions in $L^1({\mathbb B
}^{d+1};(1-|x|^2)^{\gamma_{\kappa}-1})$, and moreover, the formula
\eqref{eq5} is true for these functions too and $V_{\kappa}f\in L^1({\mathbb
S}^d;h_{\kappa}^2)$. In particular, if $\phi\in
L^1([-1,1];w_{\lambda_{\kappa}})$, then for each $y\in{\mathbb
S}^d$,
\begin{gather}\label{eq6}
\int_{{\mathbb{B}}^{d+1}}\phi(\langle
x,y\rangle)\big(1-|x|^2\big)^{\gamma_{\kappa}-1}dx
=\frac{\pi^{(d+1)/2}\Gamma(\gamma_{\kappa})}{\Gamma(\lambda_{\kappa}+1)}
\int_{-1}^1\phi(t)w_{\lambda_{\kappa}}(t)dt,
\end{gather}
i.e.\ $f(x)=\phi(\langle x,y\rangle)\in L^1({\mathbb B
}^{d+1};(1-|x|^2)^{\gamma_{\kappa}-1})$, and hence
$V_{\kappa}[\phi(\langle\cdot,y\rangle)]$ is well def\/ined and in
$L^1({\mathbb S}^d;h_{\kappa}^2)$. In addition, we have the
following symmetric relation
\begin{gather}\label{eq7}
V_{\kappa}[\phi(\langle\cdot,y\rangle)](x)=V_{\kappa}[\phi(\langle
x,\cdot\rangle)](y), \qquad \hbox{for a.e.} \quad (x,y)\in{\mathbb
S}^d\times{\mathbb S}^d.
\end{gather}
The validity of \eqref{eq7} for polynomials and for all $(x,y)\in{\mathbb
S}^d\times{\mathbb S}^d$ follows from Gegenbauer expansions and the
symmetry of the reproducing kernels $P_n(h_\kappa^2;x,y)$, If
$\phi\in L^1([-1,1];w_{\lambda_{\kappa}})$ and $\phi_1$ is a~univariate polynomial, then applying \eqref{eq5} and \eqref{eq6},
\begin{gather*}
 \int_{{\mathbb{S}}^{d}}\int_{{\mathbb{S}}^{d}}|V_{\kappa}[\phi(\langle\cdot,y\rangle)](x)-V_{\kappa}[\phi(\langle
x,\cdot\rangle)](y)|h_\kappa^2(x)h_\kappa^2(y) d\omega_d(x)d\omega_d(y)\\
\qquad{}=\int_{{\mathbb{S}}^{d}}\int_{{\mathbb{S}}^{d}}
|V_{\kappa}[(\phi-\phi_1)(\langle\cdot,y\rangle)](x)-V_{\kappa}[(\phi-\phi_1)(\langle
x,\cdot\rangle)](y)|h_\kappa^2(x)h_\kappa^2(y) d\omega_d(x)d\omega_d(y)\\
\qquad{} \le 2c_{\kappa}^{-2}\int_{-1}^1|\phi(t)-\phi_1(t)|w_{\lambda_{\kappa}}(t)dt,
\end{gather*}
which implies \eqref{eq7} by the density of polynomials in
$L^1([-1,1];w_{\lambda_{\kappa}})$. Following the above remarks, the
Funk--Hecke formula for $h$-harmonics proved in~\cite{Xu2} (for
continuous functions there only) holds also for $\phi\in
L^1([-1,1];w_{\lambda_{\kappa}})$, that is
\begin{gather}\label{eq8}
c_{\kappa}\int_{{\mathbb{S}}^{d}}
V_{\kappa}[\phi(\langle\cdot,y\rangle)](x)H_n(x)h_\kappa^2(x)
d\omega_d(x)=L_n(\phi)H_n(y)
\end{gather}
for each $H_n\in{\mathcal{H}}_{n}^{h,d+1}$ and $y\in{\mathbb S}^d$,
where
\begin{gather}\label{eq9}
L_n(\phi)=
\int_{-1}^1\phi(t)\frac{C_n^{\lambda_{\kappa}}(t)}{C_n^{\lambda_{\kappa}}(1)}w_{\lambda_{\kappa}}(t)dt.
\end{gather}

The convolution $f\ast_{\kappa}\phi$ of two functions $f\in
L^1({\mathbb S}^d;h_{\kappa}^2)$ and $\phi\in
L^1([-1,1];w_{\lambda_{\kappa}})$ is def\/ined in~\cite{Xu4}, by
\begin{gather}\label{eq10}
f\ast_{\kappa}\phi(x)=c_{\kappa}\int_{{\mathbb{S}}^{d}}f(y)
V_{\kappa}[\phi(\langle x,\cdot\rangle)](y)h_\kappa^2(y)
d\omega_d(y).
\end{gather}
The Young inequality concerning such convolution is proved in~\cite{Xu4}, that is, for $p,q,r\ge1$ with $r^{-1}=p^{-1}+q^{-1}-1$,
\begin{gather}\label{eq11}
\|f\ast_{\kappa}\phi\|_{\kappa,r}\le\|f\|_{\kappa,p}\|\phi\|_{\lambda_{\kappa},q}.
\end{gather}

A typical example of Dunkl's theory is the case when $G=Z_2^{d+1}$,
for which, the func\-tion~$h_\kappa(x)$ has the form $h_\kappa(x) =
|x_1|^{\kappa_1}\cdots |x_{d+1}|^{\kappa_{d+1}}$ and the
intertwining operator $V_{\kappa}$ is given by
\begin{gather}\label{eq12}
 V_{\kappa}f(x)=\tilde{c}_{\kappa}\int_{[-1,1]^{d+1}}
f(x_1t_1,\dots,x_{d+1}t_{d+1}) \prod_{i=1}^{d+1}(1+t_i)
(1-t_i^2)^{\kappa_i -1} dt_1\cdots dt_{d+1},
\end{gather}
where $\tilde{c}_\kappa= \tilde{c}_{\kappa_1} \cdots
\tilde{c}_{\kappa_{d+1}}$.

\section[The spherical Radon-Dunkl transform]{The spherical Radon--Dunkl transform}\label{section3}

For $f\in L^1({\mathbb S}^d;h_{\kappa}^2)$, its generalized
spherical means $M_{\tau}^{\kappa}f(x)$ due to \cite{Xu3} is def\/ined
by the equation
\begin{gather}\label{eq13}
\int_{-1}^1
M_{\tau}^{\kappa}f(x)\phi(\tau)w_{\lambda_{\kappa}}(\tau)d\tau
=f\ast_{\kappa}\phi(x)
\end{gather}
for any $\phi$ in $L^1([-1,1];w_{\lambda_{\kappa}})$. Since, for
$\phi\in L^{\infty}({\mathbb S}^d)$,
\begin{gather*}
\left|\int_{-1}^1
M_{\tau}^{\kappa}f(x)\phi(\tau)w_{\lambda_{\kappa}}(\tau)d\tau\right|\le
\|f\|_{\kappa,1}\|\phi\|_{\infty},
\end{gather*}
it follows that, for each $x\in{\mathbb S}^d$, the function
$\psi_x(\tau)=M_{\tau}^{\kappa}f(x)\in
L^1([-1,1];w_{\lambda_{\kappa}})$. This shows that for almost all
$\tau\in[-1,1]$, $M_{\tau}^{\kappa}f$ is well def\/ined. To give
further illustration of $M_{\tau}^{\kappa}$, we introduce the space
$W_m^p({\mathbb S}^d;h_{\kappa}^2)(\subseteq L^p({\mathbb
S}^d;h_{\kappa}^2))$ of functions for $m\ge0$, such that for $f\in
W_m^p({\mathbb S}^d;h_{\kappa}^2)$, there exist some $g\in
L^p({\mathbb S}^d;h_{\kappa}^2)$ satisfying
$Y_0(h_\kappa^2;f)=Y_0(h_\kappa^2;g)$ and
$[n(n+2\lambda_{\kappa})]^{m/2}Y_n(h_\kappa^2;f)=Y_n(h_\kappa^2;g)$
for all $n=1,2,\dots$. In view of~\eqref{eq4}, we formally write
$g=(-\Delta_{h,0})^{m/2}f$. It is noted that for even $m$,
$C^{m}({\mathbb S}^d)\subseteq W_m^p({\mathbb S}^d;h_{\kappa}^2)$
($1\le p\le\infty$). The following properties of $M_{\tau}^{\kappa}$
are proved in \cite{Xu3,Xu4}.

\begin{proposition} \label{proposition1} \qquad
\begin{enumerate}\itemsep=0pt
\item[$(i)$] If $f_0(x)\equiv1$, then
$M_{\tau}^{\kappa}f_0(x)\equiv1$.

\item[$(ii)$] For each $\tau\in[-1,1]$, there is an extension of
$M_{\tau}^{\kappa}$ to $L^p({\mathbb S}^d;h_{\kappa}^2)$ $(1\le
p<\infty)$, or $C({\mathbb S}^d)$ $(p=\infty)$, such that
\begin{gather*}
 \|M_{\tau}^{\kappa}f\|_{\kappa,p}\le \|f\|_{\kappa,p}, \qquad
 \tau\in[-1,1].
\end{gather*}

\item[$(iii)$] For $f\in L^1({\mathbb S}^d;h_{\kappa}^2)$,
\begin{gather*}
 Y_n(h_\kappa^2;M_{\tau}^{\kappa}f;x)=\frac{C_n^{\lambda_{\kappa}}(\tau)}{C_n^{\lambda_{\kappa}}(1)}Y_n\big(h_\kappa^2;f;x\big),
\end{gather*}
and in particular,
$\Delta_{h, 0}(M_{\tau}^{\kappa}f)=M_{\tau}^{\kappa}(\Delta_{h,0}f)$
if $\Delta_{h, 0}f\in L^1({\mathbb S}^d;h_{\kappa}^2)$.
\end{enumerate}
\end{proposition}

\begin{proof}
Here we give an independent, but simpler proof for part $(ii)$ as
follows. For all $\phi\in L^1([-1,1]$; $w_{\lambda_{\kappa}})$ and
$g\in L^{p'}({\mathbb S}^d;h_{\kappa}^2)$, it follows from \eqref{eq13} that
\begin{gather*}
\int_{-1}^1\phi(\tau)\xi(\tau)w_{\lambda_{\kappa}}(\tau)d\tau
=c_{\kappa}\int_{{\mathbb
S}^d}(f\ast_{\kappa}\phi)(x)g(x)h_{\kappa}^2(x)d\omega_d(x),
\end{gather*}
where $\xi(\tau)=c_{\kappa}\int_{{\mathbb S}^d}
M_{\tau}^{\kappa}f(x)\cdot g(x)h_{\kappa}^2(x)d\omega_d(x)$. By
using \eqref{eq7} and \eqref{eq10}, the right-hand side above becomes
$c_{\kappa}\int_{{\mathbb
S}^d}f(y)(g\ast_{\kappa}\phi)(y)h_{\kappa}^2(y)d\omega_d(y)$, and
then, by applying H\"older's inequality and the Young inequality
\eqref{eq11}, its absolute value is dominated by
$\|f\|_{\kappa,p}\|g\|_{\kappa,p'}\|\phi\|_{\lambda_{\kappa},1}$.
This gives that $\sup_{\tau\in[-1,1]}|\xi(\tau)|\le
\|f\|_{\kappa,p}\|g\|_{\kappa,p'}$, which means that, for almost all
$\tau\in[-1,1]$, $\|M_{\tau}^{\kappa}f\|_{\kappa,p}\le
\|f\|_{\kappa,p}$. If $f$ is a polynomial, part $(iii)$ implies that
$M_{\tau}^{\kappa}f(x)$ is a continuous function of
$(\tau,x)\in[-1,1]\times{\mathbb S}^d$, so that
$\|M_{\tau}^{\kappa}f\|_{\kappa,p}\le \|f\|_{\kappa,p}$ is true for
all $\tau\in[-1,1]$ in this case. Finally, from density of the set
of polynomials, for each $\tau\in[-1,1]$, $M_{\tau}^{\kappa}$ can be
extended to all functions in $L^p({\mathbb S}^d;h_{\kappa}^2)$
($1\le p<\infty$), or $C({\mathbb S}^d)$. Following this, part $(iii)$
also holds for $f\in L^1({\mathbb S}^d;h_{\kappa}^2)$ and each
$\tau\in[-1,1]$, and moreover,
\begin{gather}\label{eq14}
 M_{\tau}^{\kappa}f\sim\sum_{n=0}^{\infty}
 \frac{C_n^{\lambda_{\kappa}}(\tau)}{C_n^{\lambda_{\kappa}}(1)}Y_n(h_\kappa^2;f;x).
\end{gather}

We note that when $f$ is even in ${\mathbb S}^d$, $M_{\tau}f$ is
even for $\tau\in(-1,1)$.
\end{proof}

 The following proposition gives a pointwise description of
$M_{\tau}^{\kappa}$ for a larger class of functions.

\begin{proposition}\label{proposition2}
For $f\in W_m^2({\mathbb S}^d;h_{\kappa}^2)$ with
$m>\lambda_{\kappa}+1$, $M_{\tau}^{\kappa}f(x)$ is a continuous
function of $(\tau,x)\in[-1,1]\times{\mathbb S}^d$ and
\begin{gather*}
 M_{\tau}^{\kappa}f(x)=\sum_{n=0}^{\infty}
 \frac{C_n^{\lambda_{\kappa}}(\tau)}{C_n^{\lambda_{\kappa}}(1)}Y_n\big(h_\kappa^2;f;x\big),
\end{gather*}
the series on the right-hand side being absolutely and uniformly
convergent.
\end{proposition}

\begin{proof}
It is noted that for $t\in[-1,1]$, $|C_n^{\lambda_{\kappa}}(t)|\le
C_n^{\lambda_{\kappa}}(1)=(2\lambda_{\kappa})_n/n!\simeq
n^{2\lambda_{\kappa}-1}$ \cite[p.~19]{DuX}. From~\eqref{eq2}, for $g\in
L^2({\mathbb S}^d;h_{\kappa}^2)$ and all $x\in{\mathbb S}^d$ we have
$|Y_n(h_\kappa^2;g;x)|\le\|g\|_{\kappa,2}\|P_n(h_{\kappa};x,\cdot)\|_{\kappa,2}$.
In view of orthogonality of $h$-harmonics and from~\eqref{eq3},
\begin{gather*}
\|P_n(h_{\kappa};x,\cdot)\|_{\kappa,2}^2=
P_n(h_{\kappa};x,x)\le\lambda_{\kappa}^{-1}(n+\lambda_{\kappa})C_n^{\lambda_{\kappa}}(1)
\simeq \lambda_{\kappa}^{-1}n^{2\lambda_{\kappa}}.
\end{gather*}
Therefore $|Y_n(h_\kappa^2;g;x)|\le
c\|g\|_{\kappa,2}n^{\lambda_{\kappa}}$. If $f\in W_m^2({\mathbb
S}^d;h_{\kappa}^2)$, then for $n\ge1$,
$Y_n(h_\kappa^2;f;x)=[n(n+2\lambda_{\kappa})]^{-m/2}Y_n(h_\kappa^2;(-\Delta_{h,0})^{m/2}f;x)$,
so that $|Y_n(h_\kappa^2;f;x)|\le
c\|(-\Delta_{h,0})^{m/2}f\|_{\kappa,2}n^{\lambda_{\kappa}-m}$.
Hence, when $m>\lambda_{\kappa}+1$, the series in \eqref{eq14} converges
absolutely and uniformly for $(\tau,x)\in[-1,1]\times{\mathbb S}^d$.
In view of the uniqueness of $h$-harmonic expansion following from
its Ces\`aro summability (see~\cite{Xu1}), the conclusions in the
proposition are proved.
\end{proof}

Now we def\/ine the transform $R_{\kappa}$ by
\begin{gather*}
R_{\kappa}f=M_{0}^{\kappa}f,
\end{gather*}
and call $R_{\kappa}$ the spherical Radon--Dunkl transform. By
Proposition~\ref{proposition1}~$(ii)$, $R_{\kappa}f$ is well def\/ined for $f\in
L^1({\mathbb S}^d;h_{\kappa}^2)$, and moreover, from Propositions~\ref{proposition1}
and~\ref{proposition2}, we have the following corollary.

\begin{corollary} \label{corollary1} \qquad

\begin{enumerate}\itemsep=0pt

\item[$(i)$] For $f\in L^p({\mathbb S}^d;h_{\kappa}^2)$ $(1\le p<\infty)$,
or $f\in C({\mathbb S}^d)$ $(p=\infty)$, we have
$\|R_{\kappa}f\|_{\kappa,p}\le \|f\|_{\kappa,p}$, and
\begin{gather}\label{eq15}
 R_{\kappa}f\sim\sum_{n=0}^{\infty}b_nY_n\big(h_\kappa^2;f;x\big),
\end{gather}
where
\begin{gather}\label{eq16}
b_{n}=\left\{\begin{array}{ll}(-1)^{\frac{n}{2}}\frac{\Gamma(\lambda_{\kappa}+1/2)}{\Gamma(1/2)}
\frac{\Gamma((n+1)/2)}{\Gamma(\lambda_{\kappa}+(n+1)/2)},
& \hbox{for} \ n  \ \hbox{even};\\
 0, & \hbox{for} \ n \ \hbox{odd}.\end{array}\right.
\end{gather}

\item[$(ii)$] For $f\in W_m^2({\mathbb S}^d;h_{\kappa}^2)$ with
$m>\lambda_{\kappa}+1$, $R_{\kappa}f(x)$ is a continuous function on
${\mathbb S}^d$ and
\begin{gather*}
 R_{\kappa}f(x)=\sum_{n=0}^{\infty}b_nY_n\big(h_\kappa^2;f;x\big),
\end{gather*}
where the series on the right-hand side is absolutely and uniformly
convergent.
\end{enumerate}
\end{corollary}

The numbers
$b_{n}=C_n^{\lambda_{\kappa}}(0)/C_n^{\lambda_{\kappa}}(1)$ are
computed by using 10-9(3) and 10-9(19) in~\cite{Er}.

The following is a nontrivial example of $R_{\kappa}$.
We consider the group $G=Z_2^{d+1}$, with
$\kappa=(\kappa_1,0,\dots,0)$ and $\kappa_1>0$. In this case,
$h_\kappa(x)=|x_1|^{\kappa_1}$ and the intertwining operator
$V_{\kappa}$ in \eqref{eq12} reduces to{\samepage
\begin{gather*}
 V_{\kappa}f(x)=\tilde{c}_{\kappa_1}\int_{-1}^1
f(x_{1}t,\tilde{x})(1+t) \big(1-t^2\big)^{\kappa_1-1}dt,
\end{gather*}
where $x=(x_1,\tilde{x})$ with
$\tilde{x}=(x_2,\dots,x_{d+1})\in{\mathbb R}^d$.}

We shall show that, for $f\in C({\mathbb S}^d)$ and for $x_1\neq0$,
\begin{gather}\label{eq17}
M_{\tau}^{\kappa}f(x)=
\frac{c_{\kappa}\tilde{c}_{\kappa_1}w_{\lambda_{\kappa}}^{-1}(\tau)}{|x_1|^{2\kappa_1}}
\int_{\Omega_{\tau}}f(y)|\langle
x,y\rangle-\tau|^{\kappa_1-1}|\langle
x\sigma,y\rangle-\tau|^{\kappa_1}d\omega_d(y),
\end{gather}
where $\Omega_{\tau}=\{y\in{\mathbb S}^d: \; y=(y_1,\tilde{y}) \
\hbox{with}\ |\langle\tilde{x},\tilde{y}\rangle-\tau|<|x_1y_1|\}$,
and $\sigma$ is the ref\/lection such that
$x\sigma=(-x_1,x_2,\dots,x_{d+1})$. Indeed, from the above formula
for $V_{\kappa}$,
\begin{gather*}
V_{\kappa}[\phi(\langle x,\cdot\rangle)](y)
=\tilde{c}_{\kappa_1}\int_{-1}^1
\phi(x_{1}y_{1}t+\langle\tilde{x},\tilde{y}\rangle)(1+t)
(1-t^2)^{\kappa_1-1}dt.
\end{gather*}
When $y_1\neq 0$, taking the substitution of variables
$t=(\tau-\langle\tilde{x},\tilde{y}\rangle)/(x_1y_1)$, we get
\begin{gather*}
V_{\kappa}[\phi(\langle x,\cdot\rangle)](y)
=\frac{\tilde{c}_{\kappa_1}}{|x_1y_1|^{2\kappa_1}}
\int_{\langle\tilde{x},\tilde{y}\rangle-|x_1y_1|}^{\langle\tilde{x},\tilde{y}\rangle+|x_1y_1|}
\phi(\tau)|\langle x,y\rangle-\tau|^{\kappa_1-1}| \langle
x\sigma,y\rangle-\tau|^{\kappa_1}d\tau.
\end{gather*}
Substituting this into the def\/inition \eqref{eq10} of $f\ast_{\kappa}\phi$,
we have
\begin{gather*}
f\ast_{\kappa}\phi(x)=\int_{-1}^1
A_{\tau}f(x)\phi(\tau)w_{\lambda_{\kappa}}(\tau)d\tau
\end{gather*}
for all $\phi\in L^1([-1,1];w_{\lambda_{\kappa}})$, where
$A_{\tau}f(x)$ denotes the expression on the right-hand side of~\eqref{eq17}. Then from~\eqref{eq13}, $A_{\tau}f(x)=M_{\tau}^{\kappa}f(x)$, so that~\eqref{eq17} is proved.

Taking $\tau=0$ in \eqref{eq17}, we get that, for $f\in C({\mathbb S}^d)$
and for $x_1\neq0$,
\begin{gather*}
R_{\kappa}f(x)=
\frac{c_{\kappa}\tilde{c}_{\kappa_1}\tilde{c}_{\lambda_{\kappa}+1/2}^{-1}}{|x_1|^{2\kappa_1}}
\int_{\Omega_{0}}f(y)|\langle x,y\rangle|^{\kappa_1-1}|\langle
x\sigma,y\rangle|^{\kappa_1}d\omega_d(y).
\end{gather*}

\section[Inversion formulas for $R_{\kappa}$ by means of spherical Riesz-Dunkl potentials]{Inversion formulas for $\boldsymbol{R_{\kappa}}$\\ by means of spherical Riesz--Dunkl potentials}\label{section4}

 For $f\in L^1({\mathbb S}^d;h_{\kappa}^2)$ and $\Re\alpha>0$,
$\alpha\neq1,3,5,\dots$, we def\/ine its spherical Riesz--Dunkl
poten\-tial~$I^{\alpha}_{\kappa}f$ by
\begin{gather}\label{eq18}
I^{\alpha}_{\kappa}f(x)=C_{\kappa,\alpha}\int_{{\mathbb S}^d}
f(y)V_{\kappa}(|\langle x,\cdot\rangle|^{\alpha-1})(y)
h_{\kappa}^2(y)d\omega_d(y),
\end{gather}
where
$C_{\kappa,\alpha}=\frac{\sqrt{\pi}\Gamma((1-\alpha)/2)}{\Gamma(\lambda_{\kappa}+1)\Gamma(\alpha/2)}c_{\kappa}$.
\begin{proposition}\label{proposition3}
For $\Re\alpha>0$, $\alpha\neq1,3,5,\dots$, $I_{\kappa}^{\alpha}f$
is well defined for each $f\in L^1({\mathbb S}^d;h_{\kappa}^2)$, and
moreover, we have the following statements:
\begin{enumerate}\itemsep=0pt
\item[$(i)$] for $1\le p\le\infty$, there exists a constant $c>0$,
such that for all $f\in L^p({\mathbb S}^d;h_{\kappa}^2)$,
$\|I_{\kappa}^{\alpha}f\|_{\kappa,p}\le c\|f\|_{\kappa,p}$;

\item[$(ii)$] if the $h$-harmonic expansion of a function $f\in
L^1({\mathbb S}^d;h_{\kappa}^2)$ is given by~\eqref{eq1}, then
$I^{\alpha}_{\kappa}f$ has the following expansion
\begin{gather}\label{eq19}
I^{\alpha}_{\kappa}f(x)\sim\sum_{n=0}^\infty
b_{n,\alpha}Y_n\big(h_\kappa^2;f;x\big), \qquad x\in{\mathbb S}^d,
\end{gather}
where
\begin{gather}\label{eq20}
b_{n,\alpha}=\left\{\begin{array}{ll}
\displaystyle (-1)^{\frac{n}{2}}
\frac{\Gamma((n+1-\alpha)/2)}{\Gamma(\lambda_{\kappa}+(n+1+\alpha)/2)},
&  \hbox{for} \ n  \ \hbox{even};\vspace{2mm}\\
 0, & \hbox{for} \ n \ \hbox{odd}.\end{array}\right.
\end{gather}
\end{enumerate}
\end{proposition}

The conclusions in the proposition are contained in Proposition~2.9
of~\cite{Xu4}. Here we give a~short presentation. Since
$I_{\kappa}^{\alpha}f=cf\ast_{\kappa}\phi$ with
$\phi(t)=|t|^{\alpha-1}\in L^1([-1,1];w_{\lambda_{\kappa}})$, part
$(i)$ follows from the Young inequality~\eqref{eq11} immediately. For part
$(ii)$, since  $V_{\kappa}(|\langle\cdot,y\rangle|^{\alpha-1})\in
L^1({\mathbb S}^d;h_{\kappa}^2)$ from~\eqref{eq5} and~\eqref{eq6}, we have, using
\eqref{eq2}, \eqref{eq8} and \eqref{eq9},
\begin{gather}\label{eq21}
Y_n\big(h_\kappa^2;V_{\kappa}\big(|\langle\cdot,y\rangle|^{\alpha-1}\big);x\big)
=L_n P_n\big(h_\kappa^2;x,y\big),
\end{gather}
where
\begin{gather*}
L_n=\int_{-1}^1|t|^{\alpha-1}\frac{C_n^{\lambda_{\kappa}}(t)}{C_n^{\lambda_{\kappa}}(1)}w_{\lambda_{\kappa}}(t)dt.
\end{gather*}
From~\eqref{eq2}, \eqref{eq7} and~\eqref{eq18}, one can get
\begin{gather*}
Y_n\big(h_\kappa^2;I^{\alpha}_{\kappa}f;x\big)=
C_{\kappa,\alpha}\int_{{\mathbb S}^d}
f(z)Y_n\big(h_\kappa^2;V_{\kappa}\big(|\langle\cdot,z\rangle|^{\alpha-1}\big);x\big)
h_{\kappa}^2(z)d\omega_d(z),
\end{gather*}
and then applying~\eqref{eq21},
$Y_n(h_\kappa^2;I^{\alpha}_{\kappa}f;x)=b_{n,\alpha}Y_n(h_\kappa^2;f;x)$
with $b_{n,\alpha}=c_{\kappa}^{-1}C_{\kappa,\alpha}L_n$. It is clear
that $b_{n,\alpha}=0$ for odd $n$. When $n$ is even, we use 7.311(2)
in~\cite{GR}, part~(v) in  \cite[p.~19]{DuX}, and some properties
of the gamma function, to get the stated value of $b_{n,\alpha}$.

It is easy to see that~\eqref{eq19} and~\eqref{eq20} allow us to extend the family
$\{I_{\kappa}^{\alpha}: \; \Re\alpha>0, \alpha\neq1,3,5,\dots\}$ to
a larger one, which leads to the following def\/inition. We put
$\Pi=\{\alpha\in{\mathbb C}: \; \alpha\neq1,3,5,\dots\}$.

\begin{definition} \label{definition1} Let $\alpha\in\Pi$. For $f\in
L^1({\mathbb S}^d;h_{\kappa}^2)$, we def\/ine $I_{\kappa}^{\alpha}f$
by the following $h$-harmonic expansion
\begin{gather}\label{eq22}
I^{\alpha}_{\kappa}f\sim\sum_{n=0}^\infty
b_{n,\alpha}Y_n(h_\kappa^2;f;x), \qquad x\in{\mathbb S}^d,
\end{gather}
where $b_{n,\alpha}$ is given by~\eqref{eq20}.
\end{definition}

It is clear that $I_{\kappa}^{\alpha}f$ is well def\/ined for $f\in
C^{\infty}({\mathbb S}^d)$. In general, $I_{\kappa}^{\alpha}f$ may
be a distribution on ${\mathbb S}^d$. Since $|b_{n,\alpha}|\le
cn^{-\lambda_{\kappa}-\Re\alpha}$, then $I_{\kappa}^{\alpha}f\in
L^2({\mathbb S}^d;h_{\kappa}^2)$ when $f\in L^2({\mathbb
S}^d;h_{\kappa}^2)$ and $\Re\alpha\ge-\lambda_{\kappa}$. For
$\Re\alpha<-\lambda_{\kappa}$ and $m\ge-\lambda_{\kappa}-\Re\alpha$,
since for $n\ge1$,
$Y_n(h_\kappa^2;f;x)=[n(n+2\lambda_{\kappa})]^{-m/2}Y_n(h_\kappa^2;(-\Delta_{h,0})^{m/2}f;x)$,
we also have $I_{\kappa}^{\alpha}f\in L^2({\mathbb
S}^d;h_{\kappa}^2)$  when $f\in W_m^2({\mathbb S}^d;h_{\kappa}^2)$.
We denote by $W_{m,\,e}^2({\mathbb S}^d;h_{\kappa}^2)$ the subspace
of even functions of $W_m^2({\mathbb S}^d;h_{\kappa}^2)$.

\begin{theorem}\label{theorem1}
If $\alpha$, $-2\lambda_{\kappa}-\alpha\in\Pi$ and
$m\ge\max\{0,-\lambda_{\kappa}-\Re\alpha\}$, then
$I_{\kappa}^{\alpha}$ is an isomorphism between
$W_{m,\,e}^2({\mathbb S}^d;h_{\kappa}^2)$ and
$W_{m+\lambda_{\kappa}+\Re\alpha,\,e}^2({\mathbb
S}^d;h_{\kappa}^2)$, and
\begin{gather*}
(I_{\kappa}^{\alpha})^{-1}=I_{\kappa}^{-2\lambda_{\kappa}-\alpha}.
\end{gather*}
\end{theorem}

In fact, proceeding the above process, it is not dif\/f\/icult to show
that for $f\in W_{m,\,e}^2({\mathbb S}^d;h_{\kappa}^2)$ with
$m\ge\max\{0,-\lambda_{\kappa}-\Re\alpha\}$, we have
$I_{\kappa}^{-2\lambda_{\kappa}-\alpha}I_{\kappa}^{\alpha}f=f$. For
$f\in W_{m',\,e}^2({\mathbb S}^d;h_{\kappa}^2)$ with
$m'=m+\lambda_{\kappa}+\Re\alpha$, since
$m'\ge\max\{0,-\lambda_{\kappa}-\Re\alpha'\}$ with
$\alpha'=-2\lambda_{\kappa}-\alpha$, we again have
$I_{\kappa}^{-2\lambda_{\kappa}-\alpha'}I_{\kappa}^{\alpha'}f=f$,
i.e.\ $I_{\kappa}^{\alpha}I_{\kappa}^{-2\lambda_{\kappa}-\alpha}f=f$.
Combining the two cases proves the theorem.

To go further, for $r\in{\mathbb Z}_+$ (nonnegative integers) we
def\/ine
\begin{gather*}
P_{r,\alpha}(\Delta_{h,\,0})=
 \left \{ \begin{array}{ll}
  \text{the identity operator},\quad & r=0,\vspace{1mm}\\
\displaystyle 4^{-r}\prod\limits_{j=1}^r[-\Delta_{h,\,0}+a_j],\quad & r\ge
1,
\end{array}
\right.
\end{gather*}
where $a_j=(2{\lambda_\kappa}-2r+2j+\alpha-1)(2r-2j+1-\alpha)$.

\begin{lemma}\label{lemma1}
If $\alpha\in\Pi$, and $r\in{\mathbb Z}_+$ such that
$2r-2\lambda_{\kappa}-\alpha\in\Pi$, then for even $n$ and
$Y_n\in{\mathcal H}_n^{h,\,d+1}$,
\begin{gather*}
P_{r,\alpha}(\Delta_{h,\,0})I_\kappa^{2r-2{\lambda_\kappa}-\alpha}I_\kappa^\alpha
Y_n=Y_n.
\end{gather*}
\end{lemma}

\begin{proof}
 From \eqref{eq19} and \eqref{eq20}, we have
\begin{gather*}
I_\kappa^{2r-2{\lambda_\kappa}-\alpha}I_\kappa^\alpha Y_n =
\frac{\Gamma((n+2{\lambda_\kappa}-2r+\alpha+1)/2)\Gamma((n+1-\alpha)/2)}
     {\Gamma((n+2r-\alpha+1)/2)\Gamma((n+2{\lambda_\kappa}+\alpha+1)/2)}Y_n.
\end{gather*}
Furthermore, from \eqref{eq4},
\begin{gather*}
P_{r,\alpha}(\Delta_{h, 0})Y_n=\prod_{j=1}^r
\left(\frac{n+2{\lambda_\kappa}+\alpha-1}{2}-r+j\right)\left(\frac{n+1-\alpha}{2}+r-j\right)Y_n.
\end{gather*}
Using the properties of $\Gamma$-functions, the result is obtained.
\end{proof}

The following theorem is a direct consequence of the above lemma.
\begin{theorem}\label{theorem2}
If $\alpha\in\Pi$, and $r\in{\mathbb Z}_+$ such that
$2r-2\lambda_{\kappa}-\alpha\in\Pi$ and
$r\ge\lambda_{\kappa}+\Re\alpha/2$, then for even $f\in
C^{\infty}({\mathbb S}^d)$ and $g=I_{\kappa}^{\alpha}f$, we have the
inversion formula
\begin{gather*}
f=P_{r,\alpha}(\Delta_{h,0})I_{\kappa}^{2r-2{\lambda_\kappa}-\alpha}g.
\end{gather*}
\end{theorem}

Now we turn to the inversion problem of the spherical Radon--Dunkl
transform $R_{\kappa}$. From \eqref{eq15}, \eqref{eq16}, \eqref{eq20} and \eqref{eq22}, we see that
\begin{gather}\label{eq23}
R_{\kappa}f=\pi^{-1/2}\Gamma(\lambda_{\kappa}+1/2)I_{\kappa}^0f.
\end{gather}
This consistency  can be also seen from the following equalities
\begin{gather}\label{eq24}
\frac{\Gamma(\lambda_{\kappa}+1)\Gamma(\alpha/2)}{\sqrt{\pi}\Gamma((1-\alpha))/2}I_{\kappa}^{\alpha}f(x)=
f\ast_{\kappa}\phi=\int_{-1}^1M_{\tau}^{\kappa}f(x)\phi(\tau)w_{\lambda_{\kappa}}(\tau)d\tau
\end{gather}
in view of \eqref{eq13} and \eqref{eq18}, where $\phi(t)=|t|^{\alpha-1}\in
L^1([-1,1];w_{\lambda_{\kappa}})$. Assume that $f\in W_m^2({\mathbb
S}^d;h_{\kappa}^2)$ with $m>\lambda_{\kappa}+1$. By Proposition~\ref{proposition2},
for each $x\in{\mathbb S}^d$, $M_{\tau}^{\kappa}f(x)$ is a
continuous function of $\tau\in[-1,1]$. Dividing each part of~\eqref{eq24}
by $\Gamma(\alpha/2)$ and taking limit for $\alpha\rightarrow0+$, we
regain the relation~\eqref{eq23}.

From Theorems~\ref{theorem1} and~\ref{theorem2}, we obtain the inversion formulas for the
spherical Radon--Dunkl transform~$R_{\kappa}$.

\begin{theorem}\label{theorem3}
$R_{\kappa}$ is an isomorphism between $W_{m, e}^2({\mathbb
S}^d;h_{\kappa}^2)$ and $W_{m+\lambda_{\kappa}, e}^2({\mathbb
S}^d;h_{\kappa}^2)$ with $m\ge0$, and
\begin{gather*}
R_{\kappa}^{-1}=\frac{\sqrt\pi}{\Gamma({\lambda_\kappa}+1/2)}I_{\kappa}^{-2\lambda_{\kappa}}.
\end{gather*}
\end{theorem}

\begin{theorem} \label{theorem4} If $r\in{\mathbb Z}_+$ such that $2r-2\lambda_{\kappa}\in\Pi$
and $r\ge\lambda_{\kappa}$, then for even $f\in C^{\infty}({\mathbb
S}^d)$ and $g=R_{\kappa}f$, we have the inversion formula
\begin{gather*}
f=\frac{\sqrt\pi}{\Gamma({\lambda_\kappa}+1/2)}P_{r,0}(\Delta_{h, 0})I_{\kappa}^{2r-2{\lambda_\kappa}}g.
\end{gather*}
\end{theorem}

For a special case, we have some simple inversion formulas for
$R_{\kappa}$, which are interesting generalizations of those about
the usual spherical Radon transform (see \cite{He1, He2,Ru3}).

\begin{corollary}\label{corollary2}
If $\lambda_{\kappa}$ is a positive integer, then an even $f\in
C^{\infty}({\mathbb S}^d)$ can be recovered by
\begin{gather*}
(i) \quad f=c'P_{r,0}(\Delta_{h, 0})R_{\kappa}R_{\kappa}f,
\end{gather*}
with $r=\lambda_{\kappa}$, $c'=\pi/\Gamma({\lambda_\kappa}+1/2)^2$,
and
\begin{gather*}
P_{r,0}(\Delta_{h, 0})=
 4^{-r}\prod\limits_{j=1}^r[-\Delta_{h, 0}+(2j-1)(2r-2j+1)];
 \end{gather*}
and
\begin{gather*}
(ii) \quad f=c'' P_{r,0}(\Delta_{h, 0})\left[\int_{{\mathbb S}^d}
R_{\kappa}f(y)V_{\kappa}(|\langle x,\cdot\rangle|)(y)
h_{\kappa}^2(y)d\omega_d(y)\right],
\end{gather*}
with $r=\lambda_{\kappa}+1$,
$c''=-2\pi^{3/2}c_{\kappa}/[\Gamma(\lambda_{\kappa}+1)\Gamma({\lambda_\kappa}+1/2)^2]$,
and
\begin{gather*}
P_{r,0}(\Delta_{h, 0})=
 4^{-r}\prod\limits_{j=1}^r[-\Delta_{h, 0}+(2j-3)(2r-2j+1)].
 \end{gather*}
\end{corollary}

\section[Inversion formulas for $R_{\kappa}$ by means of associated wavelets]{Inversion formulas for $\boldsymbol{R_{\kappa}}$ by means of associated wavelets}\label{section5}

In this section, we shall use, for a suitably chosen $\psi$ def\/ined
on $[0,\infty)$, the wavelet-like transform
\begin{gather}\label{eq25}
 W_\kappa f(t,x)=f\ast_{\kappa}\psi_t(x), \qquad
 \psi_t(\tau)=t^{-1}\psi(\tau/t),
\end{gather}
for $(t,x)\in(0,\infty)\times{\mathbb S}^d$, to present the inverse
of the spherical Radon--Dunkl transform $R_{\kappa}$ and itself.
Although $R_{\kappa}$ is def\/ined implicity and the intertwining
operator $V_{\kappa}$ is involved in the def\/inition of $W_{\kappa}$,
the approaches in studying the usual spherical Radon transform (see~\cite{Ru2}, for example) could be transplanted to $R_{\kappa}$.

The f\/irst lemma below reveals a relation of the spherical
Radon--Dunkl transform $R_{\kappa}$ with the one-dimensional
fractional integral, and the second gives a representation of the
successive action of $R_{\kappa}$ and $W_{\kappa}$ to a function. We
shall use a modif\/ied notation of the fractional integral~as
\begin{gather}\label{eq26}
B_{\delta}\phi(u)=\frac{2}{\Gamma(\delta)}\int_0^{\sqrt{u}}\phi(v)\big(u-v^2\big)^{\delta-1}dv,
\qquad u>0,
\end{gather}
for $\delta>0$, which will simplify some expressions.

\begin{lemma}\label{lemma2}
For even function $f\in L^1({\mathbb S}^d;h_\kappa^2)$ and $0<s<1$,
we have
\begin{gather}\label{eq27}
M_s^\kappa(R_\kappa
f)=\frac{\lambda_{\kappa}\pi^{-1}}{w_{\lambda_\kappa}(s)}
B_{\lambda_\kappa}(M_{\tau}^{\kappa}f)\big(1-s^2\big),
\end{gather}
where the action of $B_{\lambda_{\kappa}}$ to $M_{\tau}^{\kappa}f$
is associated with $\tau$-variable.
\end{lemma}

\begin{proof}
From the product formula of the Gegenbauer polynomial
$C_{2n}^{\lambda_{\kappa}}$ (see \cite[p.~203]{DuX}), we have
\begin{gather*}
\frac{C_{2n}^{\lambda_{\kappa}}(s)}{C_{2n}^{\lambda_{\kappa}}(1)}
\frac{C_{2n}^{\lambda_{\kappa}}(0)}{C_{2n}^{\lambda_{\kappa}}(1)}
=2\int_0^1\frac{C_{2n}^{\lambda_{\kappa}}(u\sqrt{1-s^2})}{C_{2n}^{\lambda_{\kappa}}(1)}
w_{\lambda_{\kappa}-1/2}(u)du.
\end{gather*}
By Proposition~\ref{proposition1} $(iii)$, the three quotients above are the
coef\/f\/icients of a member $Y_{2n}$ in ${\mathcal{H}}_{2n}^{h,d+1}$
under action of $M_{s}^{\kappa}$, $R_{\kappa}(=M_0^{\kappa})$, and
$M_{u\sqrt{1-s^2}}^\kappa$, respectively. Therefore,
\begin{gather*}
M_s^\kappa(R_\kappa Y_{2n}) =2\int_0^1
\big(M_{u\sqrt{1-s^2}}^\kappa Y_{2n}\big)
w_{\lambda_{\kappa}-1/2}(u)du.
\end{gather*}
Making substitution of variables $u=v/\sqrt{1-s^2}$, \eqref{eq27} is proved
for $Y_{2n}$. By Proposition~\ref{proposition1}~$(ii)$, both sides of \eqref{eq27} are bounded
operators in $L^1({\mathbb S}^d;h_\kappa^2)$, and hence, the
validity of~\eqref{eq27} for general even $f\in L^1({\mathbb
S}^d;h_\kappa^2)$ follows from density of the set of $h$-harmonics.
\end{proof}

\begin{lemma}\label{lemma3}
For even $f\in L^1({\mathbb S}^d;h_\kappa^2)$ and $\psi\in
L^1([0,\infty);dx)$, we have
\begin{gather}\label{eq28}
 W_{\kappa}(R_{\kappa}f)(t,x)=\frac{2\lambda_{\kappa}}
 {\pi}
 \int_0^1 M_s^{\kappa}f(x)(B_{\lambda_{\kappa}}\psi_t)\big(1-s^2\big)ds,
\end{gather}
provided the integral on the right-hand side exists with $|f|$ and
$|\psi|$ instead of $f$, $\psi$.
\end{lemma}

\begin{proof}
From \eqref{eq13}, \eqref{eq25} and \eqref{eq27},  we have
\begin{gather*}
W_{\kappa}(R_{\kappa}f)(t,x) = \int_{-1}^1\!
M_{s}^{\kappa}(R_{\kappa}f)(x)\cdot\psi_t(s)w_{\lambda_{\kappa}}(s)ds  = \frac{2\lambda_{\kappa}}{\pi}\int_0^1\!
B_{\lambda_\kappa}(M_{\tau}^{\kappa}f)\big(1-s^2\big)\cdot\psi_t(s)ds,
\end{gather*}
and then, substituting the formula for
$B_{\lambda_\kappa}(M_{\tau}^{\kappa}f)$ from~\eqref{eq26}, and making
changes of variables, we prove the equality in~\eqref{eq28}.
\end{proof}

\begin{theorem}\label{theorem5}
 Let
 \begin{gather}
\int_0^\infty s^j \psi(s) ds=0 \qquad\text{for all}\quad j=0,2,4,\dots, 2[{\lambda_\kappa}],\label{eq29}\\
 \int_1^\infty s^\beta |\psi(s)| ds<\infty \qquad \text{for some}\quad \beta>2{\lambda_\kappa}.\label{eq30}
\end{gather}
Then for even $f\in L^p({\mathbb S}^d;h_{\kappa}^2)$ $(1\le
p<\infty)$, or $C({\mathbb S}^d)$ $(p=\infty)$, we have
\begin{gather}\label{eq31}
\lim_{\epsilon\rightarrow0+}\|T_{\epsilon}f-f\|_{\kappa,p}=0,
\end{gather}
where
\begin{gather}\label{eq32}
T_{\epsilon}f(x)=\tilde{C}_{\psi}^{-1}\int_\epsilon^\infty
\frac{(W_{\kappa}g)(t,x)}{t^{2{\lambda_\kappa}+1}}dt, \qquad
\epsilon>0,
\end{gather}
with $g=R_{\kappa}f$ and
\begin{gather}\label{eq33}
 \tilde{C}_{\psi}=
  \left \{ \begin{array}{ll}
  \displaystyle  -\frac{2\Gamma(1-\lambda_{\kappa})}{\pi}\int_0^{\infty} s^{2{\lambda_\kappa}}\psi(s)ds,
                     \quad &  \text{if}\  \ {\lambda_\kappa}\,  \bar\in \, {\mathbb N},\vspace{1mm}\\
 \displaystyle
    \frac{4(-1)^{{\lambda_\kappa}+1}}{\pi\Gamma(\lambda_{\kappa})}
       \int_0^{\infty} s^{2{\lambda_\kappa}}\psi(s)\log s  ds,\quad &  \text{if}\  \ {\lambda_\kappa}\in{\mathbb N}.
\end{array}
\right.
\end{gather}
In addition, $\lim\limits_{\epsilon\rightarrow0+}T_{\epsilon}f(x)=f(x)$ for
almost all $x\in{\mathbb S}^d$.
\end{theorem}

\begin{proof}
Under the assumptions, by \cite[Lemma 4.12]{Ru1}, we have
$\int_0^{\infty}|B_{\lambda_{\kappa}}\psi(s)|ds<\infty$. To prove~\eqref{eq31} in general, we only need to show that it is valid for
$Y_{2n}\in{\mathcal{H}}_{2n}^{h,d+1}$, and
\begin{gather}\label{eq34}
\|T_{\epsilon}f\|_{\kappa,p}\le c\|f\|_{\kappa,p}, \qquad
\epsilon>0,
\end{gather}
where the constant $c$ is independent of $\epsilon$. The key step is
to rewrite $T_{\epsilon}$ into a convolution operator with an
approximate identity, that is,
\begin{gather}\label{eq35}
T_{\epsilon}f(x)=\frac{2\lambda_\kappa(\lambda_{\kappa}+1)}{\pi(2\lambda+1)\tilde{C}_{\psi}}
f\ast_{\kappa}K_{\epsilon},
\end{gather}
where
\begin{gather}\label{eq36}
K_\epsilon(\tau)=[w_{\lambda_{\kappa}+1}(\tau)]^{-1}
(B_{\lambda_{\kappa}+1}\psi)\left(\epsilon^{-2}\big(1-\tau^2\big)\right).
\end{gather}
Indeed, applying Lemma~\ref{lemma3} to~\eqref{eq32} gives that
\begin{gather}\label{eq37}
 \tilde{C}_{\psi}T_{\epsilon}f(x)=\frac{2\lambda_{\kappa}}
 {\pi} \int_0^1 M_s^{\kappa}f(x)\tilde{K}_{\epsilon}(s)ds,
\end{gather}
where $\tilde{K}_{\epsilon}(s)=\int_\epsilon^\infty
\frac{(B_{\lambda_{\kappa}}\psi_t)(1-s^2)}{t^{2{\lambda_\kappa}+1}}dt$.
Inserting the formula of $B_{\lambda_{\kappa}}\psi_t$ from~\eqref{eq26}, and
then, making changes of variables by $t=\xi^{-1/2}$ and
$v=\eta \xi^{-1/2}$, we get
\begin{gather*}
\tilde{K}_{\epsilon}(s)=\frac{1}{\Gamma(\lambda_{\kappa})}\int_0^{\epsilon^{-2}}
\int_0^{\sqrt{\xi(1-s^2)}}\psi(\eta)[\xi(1-s^2)-\eta^2]^{\lambda_{\kappa}-1}d\eta
d\xi.
\end{gather*}
Changing order of the integrals, it follows that
$\tilde{K}_{\epsilon}(s)=\frac{\lambda_{\kappa}+1}{2\lambda_{\kappa}+1}K_{\epsilon}(s)w_{\lambda_{\kappa}}(s)$.
Substituting this into~\eqref{eq37} and using~\eqref{eq13} yields~\eqref{eq35}.

By Lemma~2.4 in~\cite{Ru2}, we have
$\int_0^{\infty}u^{-1}(B_{\lambda_{\kappa}+1}\psi)(u)du=\pi
\tilde{C}_{\psi}/\lambda_{\kappa}$, and
$B_{\lambda_{\kappa}+1}\psi(u)=O(u^{\lambda_{\kappa}})$ for
$0<u\le1$, and $O(u^{-\rho})$ for $u>1$ with some $\rho>0$. From
these and in view of~\eqref{eq9} and~\eqref{eq36}, it follows that
\begin{gather}
L_{2n}(K_{\epsilon}) =
\frac{2\lambda_{\kappa}+1}{\lambda_{\kappa}+1}
\int_0^1(B_{\lambda_{\kappa}+1}\psi)\left(\frac{1-s^2}{\epsilon^2}\right)
\frac{C_{2n}^{\lambda_{\kappa}}(s)ds}{(1-s^2)C_{2n}^{\lambda_{\kappa}}(1)}\nonumber\\
 \phantom{L_{2n}(K_{\epsilon})}{} = \frac{\lambda_{\kappa}+1/2}{\lambda_{\kappa}+1}
\int_0^{\epsilon^{-2}}u^{-1}(B_{\lambda_{\kappa}+1}\psi)(u)
\frac{C_{2n}^{\lambda_{\kappa}}(\sqrt{1-\epsilon^2u})}{\sqrt{1-\epsilon^2u}C_{2n}^{\lambda_{\kappa}}(1)}du,
\label{eq38}
\end{gather}
which approaches to
$\frac{\pi(\lambda_{\kappa}+1/2)}{\lambda_{\kappa}(\lambda_{\kappa}+1)}\tilde{C}_{\psi}$
as $\epsilon\rightarrow+0$. For
$Y_{2n}\in{\mathcal{H}}_{2n}^{h,d+1}$, from ~\eqref{eq8} we have
$Y_{2n}\ast_{\kappa}K_{\epsilon}=L_{2n}(K_{\epsilon})Y_{2n}$, and by~\eqref{eq35}, $\lim\limits_{\epsilon\rightarrow0+}T_{\epsilon}Y_n=Y_n$ uniformly on
${\mathbb S}^d$.

To prove~\eqref{eq34}, by~\eqref{eq11}, it suf\/f\/ices to show
$\|K_{\epsilon}\|_{\lambda_{\kappa},1}\le c$ uniformly for
$\epsilon>0$ (essentially for $0<\epsilon\le1$). In fact, similarly
to~\eqref{eq38}, we have
$\|K_{\epsilon}\|_{\lambda_{\kappa},1}=L_0(|K_{\epsilon}|)$,
approaching to
\begin{gather*}
\frac{\lambda_{\kappa}+1/2}{\lambda_{\kappa}+1}
\int_0^{\infty}u^{-1}|(B_{\lambda_{\kappa}+1}\psi)(u)|du<+\infty,
\end{gather*}
as $\epsilon\rightarrow+0$. Thus \eqref{eq34}, and so~\eqref{eq31}, are proved.

In order to prove $T_{\epsilon}f$ to be convergent almost
everywhere, we need the associated maximal function
$T_*f(x)=\sup\limits_{0<\epsilon\le1}|T_{\epsilon}f(x)|$. We shall show
that $T_*f$ is dominated by the maximal function introduced in~\cite{Xu5}
\begin{gather*}
{\mathcal M}_\kappa f(x)=\sup_{0<\theta\le\pi} \frac{\int_0^\theta
(M_{\cos\varphi}^{\kappa}|f|)(x)(\sin\varphi)^{2\lambda_\kappa}\,d\varphi}{\int_0^\theta
(\sin\varphi)^{2\lambda_{\kappa}}\,d\varphi},
\end{gather*}
for $f\in L^1({\mathbb S}^d;h_\kappa^2)$, that is
\begin{gather}\label{eq39}
T_*f(x)\le c{\mathcal M}_\kappa f(x), \qquad x\in{\mathbb
S}^d.
\end{gather}
The pointwise estimates of $B_{\lambda_{\kappa}+1}\psi(u)$ can be
written as
$B_{\lambda_{\kappa}+1}\psi(u)=O(u^{\lambda_{\kappa}}(u+1)^{-\lambda_{\kappa}-\rho})$,
which implies the following estimate for $K_{\epsilon}(\cos\theta)$
\begin{gather*}
K_{\epsilon}(\cos\theta)=O(m_{\epsilon}(\theta)),  \qquad
m_{\epsilon}(\theta)=\frac{\epsilon^{2\rho}(\sin\theta)^{-1}}{(\epsilon+\sin\theta)^{2\lambda_{\kappa}+2\rho}},
\end{gather*}
with $\rho>0$. The function $m_{\epsilon}(\theta)$ does not suit the
process of integration by part in the proof of Theorem~2.6 in~\cite{Xu5}, since $m(0)=0$. Here we give a proof for the case.

From~\eqref{eq13} and~\eqref{eq35},
\begin{gather}\label{eq40}
|T_{\epsilon}f(x)|\le
c\int_{0}^{\pi/2}(M_{\cos\theta}^{\kappa}|f|)(x)m_{\epsilon}(\theta)(\sin\theta)^{2\lambda_\kappa} d\theta,
\end{gather}
where the evenness of $M_{\tau}^{\kappa}f$ is used. Splitting the
interval $[0,\pi/2]$ into $\bigcup_j[2^j\epsilon, 2^{j+1}\epsilon]$,
we evaluate each integral $U_j=\int_{2^j\epsilon}^{2^{j+1}\epsilon}$
separately. For $j\le0$, since $\epsilon+\sin\theta\asymp\epsilon$,
we have
\begin{gather*}
U_j \le
\frac{c2^{-j}}{\epsilon^{2\lambda_{\kappa}+1}}\int_{0}^{2^{j+1}\epsilon}(M_{\cos\theta}^{\kappa}|f|)(x)
(\sin\theta)^{2\lambda_\kappa} d\theta
 \le c2^{2\lambda_{\kappa}j}{\mathcal M}_\kappa f(x);
\end{gather*}
and for $j>0$, since $\epsilon+\sin\theta\asymp\theta$,
\begin{gather*}
U_j \le
\frac{c\epsilon^{2\rho}}{(2^j\epsilon)^{2\lambda_{\kappa}+2\rho+1}}
\int_{0}^{2^{j+1}\epsilon}(M_{\cos\theta}^{\kappa}|f|)(x)
(\sin\theta)^{2\lambda_\kappa} d\theta
 \le c2^{-2\rho j}{\mathcal M}_\kappa f(x).
\end{gather*}
Collecting these estimates into~\eqref{eq40} yields~\eqref{eq39}.

By Theorem~2.1 in~\cite{DX}, $T_*$ is of weak~(1,1), and strong
$(p,p)$ boundedness. Combining with the uniformly convergence of
$T_{\epsilon}$ for $h$-harmonics, for general $f\in L^1({\mathbb
S}^d;h_\kappa^2)$, $T_{\epsilon}f$ converges to $f$ almost
everywhere. The proof of Theorem~\ref{theorem5} is completed.
\end{proof}

In the following, we state two theorems, without proof, which are
analogs of Theorems~1.2 and~1.4 in~\cite{Ru2}. One is about the
reproducing property of the spherical Radon--Dunkl transform~$R_{\kappa}$, and the other illustrates the range
$R_{\kappa}(L^1({\mathbb S}^d;h_\kappa^2))$.

\begin{theorem}\label{theorem6}
Let
\begin{gather*}%\label{8}
\int_0^\infty \psi(s) ds=0,\qquad  \int_0^\infty|\psi(s)\log
s| ds<\infty.
\end{gather*}
Then for $f\in L^p({\mathbb S}^d;h_{\kappa}^2)$ $(1\le p<\infty)$,
or $C({\mathbb S}^d)$ $(p=\infty)$, we have
\begin{gather*}
\lim_{\epsilon\rightarrow0+}\|\tilde{T}_{\epsilon}f-R_{\kappa}f\|_{\kappa,p}=0,
\end{gather*}
where
$\tilde{T}_{\epsilon}f(x)=\bar{C}_{\psi}^{-1}\int_\epsilon^\infty
t^{-1}(W_{\kappa}f)(t,x)dt$  $(\epsilon>0)$, with
$\bar{C}_{\psi}=2c_{\lambda_\kappa}\int_0^\infty \psi(s)\log\frac 1
s ds$.
\end{theorem}

\begin{theorem}\label{theorem7}
 Let $\psi$ satisfy conditions~\eqref{eq29} and~\eqref{eq30},
 $g\in L^p({\mathbb S}^d;h_{\kappa}^2)$ $(1\le p<\infty)$, or
$C({\mathbb S}^d)$ $(p=\infty)$, and $\tilde{C}_{\phi}\neq 0$ be
 the constant in~\eqref{eq33}. Then the following statements are equivalent:
\begin{enumerate}\itemsep=0pt
\item[$(i)$] $g\in R_\kappa(L^p({\mathbb S}^d;h_\kappa^2))$;

\item[$(ii)$] the integrals
 $S_{\epsilon}g=\int_\epsilon^\infty t^{-2{\lambda_\kappa}-1}
(W_\kappa g)(t,x) dt
 $
 converge in the $L^p({\mathbb S}^d;h_\kappa^2)$-norm.
\end{enumerate}

If $1<p<\infty$, then $(i)$ and $(ii)$ are equivalent to
\begin{enumerate}\itemsep=0pt
\item[$(iii)$] $\sup\limits_{\epsilon>0}\|S_\epsilon
g\|_{\kappa, p}<\infty$.
\end{enumerate}
\end{theorem}

\subsection*{Acknowledgments}
This work is supported by the National Natural
 Science Foundation of China (No.~10571122),
the Beijing Natural Science Foundation, the Project of Excellent
Young Teachers and the Doctoral Programme Foundation of National
Education Ministry of China, and the Project of Beijing Education
Ministry.

\pdfbookmark[1]{References}{ref}
\LastPageEnding

\end{document}